# VARIATIONAL MONTE CARLO — BRIDGING CONCEPTS OF MACHINE LEARNING AND HIGH DIMENSIONAL PARTIAL DIFFERENTIAL EQUATIONS

MARTIN EIGEL, REINHOLD SCHNEIDER, PHILIPP TRUNSCHKE, SEBASTIAN WOLF

ABSTRACT. A statistical learning approach for parametric PDEs related to Uncertainty Quantification is derived. The method is based on the minimization of an empirical risk on a selected model class and it is shown to be applicable to a broad range of problems. A general unified convergence analysis is derived, which takes into account the approximation and the statistical errors. By this, a combination of theoretical results from numerical analysis and statistics is obtained. Numerical experiments illustrate the performance of the method with the model class of hierarchical tensors.

## 1. INTRODUCTION

In this article we explore connections between two very active research areas, namely *Machine Learning* (ML) and *Uncertainty Quantification* (UQ) *with high-dimensional partial differential equations* (PDEs). A central goal is to illustrate how ideas from statistical learning can be exploited for the solution of high-dimensional and possibly non-linear parametric deterministic problems. We employ *empirical risk minimization* (ERM) to approximate the solution of these PDEs in non-linear model classes. Although we use hierarchical tensor networks (as described e.g. in [1]) in our experiments, the presented approach is sufficiently general to also cover more complicated model classes such as deep neural networks [2]. Moreover, the framework can be applied to a large variety of problems, not necessarily related to simulation problems and differential equations. In fact, most of the previous research in this area is concerned with statistical problems which are founded on a fixed set of measured data points with statistical errors. This is fundamentally different from the measurements we consider in this work since we assume that arbitrary samples can be generated by simulation

*Date*: October 2, 2018.

Research of S. Wolf was funded in part by the DFG MATHEON project SE10.

of a computable model which – apart from numerical approximation errors – do not exhibit any statistical errors.

The term *Variational Monte Carlo* method, which we use for the presented approach, has its origins in ground state computation in Quantum Physics. In large parts (in particular for linear models), the notion coincides with our method [3].

The application focus of this work is a class of computationally demanding problems, which play an important role in UQ [4, 5]. There, the uncertainty of model data is commonly described by random fields, leading to high-dimensional parametrized PDEs. These may either be solved by Monte Carlo sampling, allowing for the computation of *quantities of interest* (i.e. functionals) of the solution, or by spectral methods leading to a *functional approximation* of the solution [6, 7]. The usually faster convergence and more detailed solution information gained by the latter approach comes at the cost of a much higher computational complexity, in particular when using a stochastic (Galerkin) projection [7–10]. With the present work, we strive to combine function space approximations (in our case in a hierarchical tensor representation) and efficient solution sampling by means of a *learning method*, i.e., the functional solution is learned from generated sampling data. The solution approximation $\Phi$ is determined by minimizing an objective functional $\mathcal{J}(\Phi) = \int_\Omega \ell(\Phi; x) \, \mathrm{d}\rho(x)$ subject to some loss function $\ell$, which e.g. yields a least squares optimization. Since $\mathcal{J}$ cannot be determined exactly in many applications such as non-linear problems, one has to resort to an empirical minimization where the computation is based on (quasi) random samples constituting the empirical risk $\mathcal{J}_N(\Phi) = \frac{1}{N} \sum_{i=1,\ldots,N} \ell(\Phi; x^i)$.

As with other numerical techniques, a convergence analysis and an a posteriori control of the occurring errors is of crucial interest. These errors can be split into a deterministic and a stochastic part, which makes it possible to combine results from numerical analysis and probability theory. Since empirical (i.e. random) quantities are involved, one does in general not obtain worst case error estimates for the stochastic error. We hence pursue the concept of *convergence in probability*. This results in error estimates which hold either with high probability or with a certain confidence, which can be improved exponentially fast by increasing the number of samples. The underlying theory has been developed in statistics and machine learning for regression and classification problems [11] and carries over with slight modifications. However, in addition to the mentioned qualitative difference of the used samples, while in statistics the probability distribution of the data is unknown and its determination can be seen as part of the task, in our case it is known and may be exploited in the method [12–15]. Another distinction is that we are primarily interested in the accuracy of the approximate solution of the PDE problem rather than computing the minimal risk.

When solving PDEs, the problem can often be cast into minimizing a specific integral like the error, the residual or the Dirichlet energy. For high-dimensional problems, due to the exponentially growing computational complexity (known as the *curse of dimensionality*), these integrals only become feasible with an empirical approach. In many situations our method is equivalent to the solution of a regression problem and we are not the first to apply regression techniques for solving PDEs, see e.g. [16]. Nevertheless, it is our intention to present a unified and general theoretical foundation in the chosen framework.

It should be noted that the considered parametric PDEs are well understood in terms of regularity and sparsity [17–21]. As such, in our opinion they represent a valuable and fruitful class of benchmark problems for machine learning algorithms, in particular since numerical methods from UQ are available as a reference.

Fundamental convergence results in statistical learning theory have been developed starting with the pioneering work of Vapnik and Chervonenkis [22–24] and Vailant [25]. We also refer to the works of Hauser [26] and Bartlett [27]. A modern treatment can be found in the recent monograph [28]. This theory is motivated by the binary classification problem and the inherent complexity is typically measured by the Vapnik-Chervonenkis (VC) dimension.

In contrast to that, a statistical learning theory where the complexity is measured in terms of the *covering number* is developed in Cucker and Smale [29] and Cucker and Zhou [30]. Since the covering number is a fundamental complexity measure in approximation theory [31], we find this approach more amenable for our approach and hence pursue this notion. Another view on regression problems in this framework can be found in [32]. While following the treatment of [29, 30], we also consider non-convex model classes.

As discussed at the end of this paper, the theoretical results so far are not optimal, at least for linear models, which is illustrated by the numerical experiments in Section 3. The derived error bounds seem to be too pessimistic and certainly can be improved for a wide class of models. A recent example of this is provided in [13] and is briefly discussed in the outlook.

The structure of the paper is as follows: In the next section, we introduce the general framework of empirical risk minimization. In Section 3, we discuss typical choices of model classes and cost functionals. Section 4 is used to derive error bounds for deterministic and stochastic error components. In Section 5 we present an application of the theory in the context of uncertainty quantification, in particular with high-dimensional PDEs. Numerical experiments illustrate the performance of the approach in Section 6. Finally, we assess our results and point in directions for future research in Section 7.

## 2. Variational Setting

Let $\mathcal{V}$ be a Hilbert space and $(\Omega, \Sigma, \rho)$ a measure space with finite measure $\rho$. We assume a given *loss function* $\ell : \mathcal{V} \times \Omega \to \mathbb{R}$ such that $\ell(\bullet\,; x)$ is continuous for all $x \in \Omega$ and $\ell(\Phi; \bullet)$ is integrable with respect to the measure $\rho$ for every $\Phi \in \mathcal{V}$. Note that we use a modified version of the usual definition of the loss function from machine learning. The objective is to find a minimizer

$$\Phi^* \in \arg\min_{\Phi \in \mathcal{V}} \mathcal{J}(\Phi)$$

of the *cost functional*

$$\mathcal{J}(\Phi) := \int_\Omega \ell(\Phi; x) \,\mathrm{d}\rho(x) \,.$$

In applications in numerical analysis the space $\mathcal{V}$ is often infinite dimensional and one hence has to confine the minimization to a discrete *model or hypothesis class* $\mathcal{M} \subseteq \mathcal{V}$. To ensure the existence of the optimum, we further assume compactness of $\mathcal{M}$. The solution of the corresponding minimization problem is denoted by

$$\Phi^{*(\mathcal{M})} \in \arg\min_{\Phi \in \mathcal{M}} \mathcal{J}(\Phi) \,.$$

Since $\rho$ is a finite measure it can be interpreted as a scaled probability measure. Accordingly, we rephrase the cost functional with a scaled version of $\ell$ as

$$\mathcal{J}(\Phi) = \int_\Omega \ell(\Phi; x) \,\mathrm{d}\rho(x) = \mathbb{E}[\ell(\Phi; \bullet)] \,.$$

We assume that independent samples $\{x^i\}_{i \leq N}$, distributed according to $\rho$, can be generated. Instead of computing $\mathcal{J}$, one can then resort to computing the *empirical cost functional*

$$\mathcal{J}_N(\Phi) := \mathbb{E}[\ell(\Phi; \bullet); N] = \frac{1}{N} \sum_{i=1}^N \ell(\Phi; x^i) \,.$$

This Monte Carlo integral provides a surrogate functional that can be minimized to obtain

$$\Phi^{*(\mathcal{M},N)} \in \underset{\Phi \in \mathcal{M}}{\arg\min}\, \mathcal{J}_N(\Phi)\ .$$

We henceforth assume that $\Phi^*$ indeed exists. The existence of the other two minimizers is guaranteed since $\mathcal{M}$ is compact and $\ell$ is continuous in its first argument.

The central topic of this paper is to analyse the errors that are introduced by restricting the optimization from $\mathcal{V}$ to $\mathcal{M}$ and by substituting $\mathcal{J}$ by the Monte Carlo surrogate $\mathcal{J}_N$. We aim to present examples that illustrate the versatility of this approach and to provide a rigorous error analysis for this general framework. The algorithm that is used to solve the minimization problem numerically can be chosen freely and we consider neither its error analysis nor its complexity.

## 3. Cost Functionals and Model Classes

This section is devoted to examples for cost functionals and model classes in our setting. For this we choose $\mathcal{V}$ as a subspace of $L^2(D)$ on a given Lipschitz domain $D$ and consider the elliptic operator equation

$$L\Phi^* = f$$

with a $\mathcal{V}$-elliptic bounded linear operator $L : \mathcal{V} \to \mathcal{V}^*$ and $f \in \mathcal{V}^*$. Common choices for $L$ are second order differential operators like the Laplace operator, the identity operator $L\Phi(x) = \Phi(x)$ or a multiplication operator $L\Phi(x) = a(x)\Phi(x)$ for some $a : D \to \mathbb{R}$.

Note that $\Phi^*$ indeed exists and that since $L$ is elliptic, an explicit restriction $\mathcal{M} \subseteq \mathcal{B}_R(0) \coloneqq \{\Phi \in \mathcal{V} \mid \|\Phi\| \leq R\}$ is not necessarily required. An upper bound $R$ for $\|\Phi^*\|_\mathcal{V}$ is given by the Lax-Milgram theorem. Depending on the model class $\mathcal{M}$ and the cost functional $\mathcal{J}$, it may however be that $\|\Phi^{*(\mathcal{M})}\|_\mathcal{V} \geq \|\Phi^*\|_\mathcal{V}$. Therefore, an additional restriction can be applied and can then be interpreted as a regularization of the problem.

Note that even though this example is a linear problem, the presented approach also allows to solve non-linear equations.

3.1. **Model classes.** The model class $\mathcal{M}$ for the optimization can be chosen quite liberally and some examples are given in the following.
  (i) Let $\mathcal{M} = \mathcal{V}$ be a finite element space.
 (ii) Let $\mathcal{M} = \mathcal{V}$ be the *Reproducing Kernel Hilbert Space (RKHS)* corresponding to a kernel $k \colon \mathcal{X}^2 \to \mathbb{R}$. Assume that the loss function has the form $\ell(\Phi; x) = \widetilde{\ell}(x, \Phi(x)) + g(\|\Phi\|_\mathcal{V})$ for a strictly increasing function $g\colon \mathbb{R} \to \mathbb{R}$. This implies that $\mathcal{J}(\Phi) = \mathbb{E}[\widetilde{\ell}(\bullet, \Phi(\bullet))] + g(\|\Phi\|)$. Then, by the *representer theorem*, $\Phi^{*(\mathcal{M},N)} = \Phi^{*(\mathcal{V},N)} \in \mathrm{span}\{k(\bullet, x_i)\}_{i \le N}$. This means that the optimization can be carried out in finite dimensions even though the solution space is indeed infinite. This is a popular choice in machine learning where the additional summand $g(\|\bullet\|_\mathcal{V})$ is often used to replace a constraint on the model class $\mathcal{M} \subseteq \{\Phi \in \mathcal{V}\colon \|\Phi\|_\mathcal{V} \le R\}$. We refer to [33] for an exposition of this topic.

Both model classes are linear but the presented framework also allows for non-linear parametrizations of the minimizer.

(iii) $\mathcal{M}$ may be a polynomial or a Gaussian mixture model.
(iv) $\mathcal{M}$ may be a set of tensors in a given tensor format. This results in multilinear models which are built upon finite dimensional subspaces of univariate functions. We point out that recently introduced hierarchical tensor representations (HT and TT formats) exhibit striking mathematical properties. For instance, they form algebraic varieties, see e.g. [1], and have an exponential power of expressiveness [34].
(v) $\mathcal{M}$ can be chosen to be multi-layer neural networks (NNs). Recent progress in machine learning demonstrates the superiority of deep neural networks over classical shallow architectures. Their theory however is rather incomplete as yet.

3.2. **Cost functionals.** The framework introduced in Section 2 is quite generic in the choice of the cost functional $\mathcal{J}$. Some possible choices are listed in the following.

(a) A *least squares approximation* of $\Phi^* = L^{-1}f$ in $\mathcal{M}$ is computed via the cost functional
$$\mathcal{J}(\Phi) = \int_\Omega |\Phi(x) - \Phi^*(x)|^2 \, \mathrm{d}\rho(x) \ .$$
The respective loss is defined by $\ell(\Phi; x) = |\Phi(x) - \Psi(x)|^2$ and the probability space is $(D, \mathcal{B}(D), \rho)$ for any measure $\rho$ that is absolutely continuous w.r.t. the Lebesgue measure. The empirical problem is equivalent to the statistical least squares regression problem with given data points $(x_i, \Phi^*(x_i))$.

(b) The *minimal residual problem* in $L^2(\Omega)$ is solved by
$$\mathcal{J}(\Phi) = \int_\Omega |L\Phi(x) - f(x)|^2 \rho(x) \, \mathrm{d}x \ .$$

The respective loss is defined by $\ell(\Phi; x) = |L\Phi(x) - f(x)|^2$ and the probability space can be chosen as above.

(c) If $L$ is self-adjoint w.r.t. the $L^2(\Omega, \rho)$ inner product, we can minimize the *Dirichlet energy*

$$\mathcal{J}(\Phi) := \frac{1}{2}(L\Phi, \Phi)_{L^2(\Omega,\rho)} - (f, \Phi)_{L^2(\Omega,\rho)}$$
$$= \frac{1}{2}\int_\Omega L\Phi(x)\Phi(x)\rho(x)\,\mathrm{d}x - \int_\Omega f(x)\Phi(x)\rho(x)\,\mathrm{d}x,$$

with $\ell(\Phi; x) := \frac{1}{2}L\Phi(x)\Phi(x) - f(x)\Phi(x)$ and the probability space again chosen as above. If $L$ is positive definite as well, then there exists an operator $B = L^{1/2}$ such that $L$ can be written as $L = B^*B$ and we can also consider

$$\mathcal{J}(\Phi) = \frac{1}{2}\int_\Omega B\Phi(x)B\Phi(x)\rho(x)\,\mathrm{d}x - \int_\Omega f(x)\Phi(x)\rho(x)\,\mathrm{d}x,$$

where $\ell(\Phi, x) := \frac{1}{2}B\Phi(x)B\Phi(x) - f(x)\Phi(x)$. A common example for this is given by $B := \sqrt{\kappa(x)}\nabla_x$. Note that although both cost functionals are equal, the corresponding empirical functionals are quite different.

**Remark 3.1.** *For many examples the exact minimum $\Phi^*$ is attained if $0 \leq \ell(\Phi^*, x) = 0$ for almost all $x$. In this case the choice of $\rho$ does not matter and we can replace $\rho$ by any density $\tilde{\rho} \ll \rho$ that is absolutely continuous w.r.t. $\rho$. The minimizers on the restricted model class $\mathcal{M} \subseteq \mathcal{V}$ can however be different from the exact minimizers. In these cases we have to keep in mind that the choice of $\tilde{\rho}$ will influence the solution $\Phi^{*(\mathcal{M})}$. This can, on the other hand, be brought to bear by weighting the samples in the empirical functional and thereby reducing instabilities.*

The preceding list is of course incomplete and is only intended to illustrate the use of surrogate functionals. Other examples like classification with softmax parametrization may prove to be interesting but are deferred to a forthcoming work.

## 4. Convergence Analysis

This central section examines the convergence of the depicted framework in terms of the errors

$$\mathcal{E}_{\text{cost}} := |\mathcal{J}(\Phi^*) - \mathcal{J}(\Phi^{*(\mathcal{M},N)})| \quad \text{and} \quad \mathcal{E}_{\text{norm}} := \|\Phi^* - \Phi^{*(\mathcal{M},N)}\|_\mathcal{V}.$$

Here, $\mathcal{E}_{\text{cost}}$ measures the error of the empirical approximation $\Phi^{*(\mathcal{M},N)}$ with respect to the cost functional $\mathcal{J}$ while $\mathcal{E}_{\text{norm}}$ determines the approximation error of the minimizer in terms of a norm $\|\bullet\|_\mathcal{V}$ related to the problem. For the further analysis, $\mathcal{E}_{\text{cost}}$ is decomposed into two parts,

$$\begin{aligned}\mathcal{E}_{\text{cost}} &\leq \mathcal{E}_{\text{appr}} + \mathcal{E}_{\text{gen}} \\ &= |\mathcal{J}(\Phi^*) - \mathcal{J}(\Phi^{*(\mathcal{M})})| \\ &\quad + |\mathcal{J}(\Phi^{*(\mathcal{M})}) - \mathcal{J}(\Phi^{*(\mathcal{M},N)})| \ .\end{aligned}$$

The first term $\mathcal{E}_{\text{appr}}$ is called the *approximation error*. It is a purely deterministic quantity due to the choice of the model class. The second term $\mathcal{E}_{\text{gen}}$ is called the *generalization error* and is a result of the use of the empirical loss functional in the minimization. This splitting allows to use of the best bounds available for each part of the error and the problem at hand. In the following, we provide generic bounds for both $\mathcal{E}_{\text{appr}}$ and $\mathcal{E}_{\text{gen}}$. In numerical applications the error of the cost functional $\mathcal{E}$ is of minor importance and one is instead interested primarily in an upper bound for the *norm error* $\mathcal{E}_{\text{norm}}$. We derive such a bound in Section 4.3.

In the subsequent analysis, certain assumptions are required, which are introduced upfront.

- **Boundedness**: There exists $C_1 > 0$ s.t. for all $\Phi \in \mathcal{M}$ it holds

$$|\ell(\Phi; x)| \leq C_1 \quad \text{for almost all } x \in \Omega \ . \tag{A1}$$

- **Lipschitz continuity on $\mathcal{M}$**: There exists $C_2 > 0$ s.t. for all $\Phi_1, \Phi_2 \in \mathcal{M}$ it holds

$$|\ell(\Phi_1; x) - \ell(\Phi_2; x)| \leq C_2 \|\Phi_1 - \Phi_2\|_{\mathcal{V}} \quad \text{for almost all } x \in \Omega \ . \tag{A2}$$

- **Global Lipschitz continuity**: There exists $C_2 > 0$ s.t. for all $\Phi_1, \Phi_2 \in \mathcal{V}$ it holds

$$|\mathcal{J}(\Phi_1) - \mathcal{J}(\Phi_2)| \leq C_2 \|\Phi_1 - \Phi_2\|_{\mathcal{V}} \ . \tag{A3'}$$

- **Bounded second derivative**: $\mathcal{J}$ is twice differentiable with bounded second derivatives, i.e.,

$$\Gamma = \sup_{\xi \in \mathcal{V}} \|\mathrm{D}^2 \mathcal{J}(\xi)\|_{\mathcal{L}(\mathcal{V},\mathcal{V}^*)} < \infty \ . \tag{A3''}$$

- **Local strong convexity**: $\mathcal{J}$ is strongly convex in a neighbourhood $\mathcal{U}$ of $\Phi^*$. This means that for all $\Phi, \Psi \in \mathcal{U}$ it holds

$$\mathcal{J}(\Phi) \geq \mathcal{J}(\Psi) + \mathrm{D}\mathcal{J}(\Psi)(\Phi - \Psi) + \frac{\gamma}{2}\|\Phi - \Psi\|_{\mathcal{V}}^2 \ . \tag{A4}$$

Assumptions (A1) and (A2) are required for bounding the generalization error in Section 4.2. Either assumption (A3') or (A3") is needed to provide bounds for the approximation error examined in Section 4.1. The last assumption (A4) is employed to bound the norm error $\mathcal{E}_{\text{norm}}$.

4.1. **Approximation error.** Define the best approximation error of $\Phi^*$ in $\mathcal{M}$ by

$$\mathcal{E}_{\text{best}} := \inf_{\Phi \in \mathcal{M}} \|\Phi^* - \Phi\|_{\mathcal{V}} \, .$$

We can bound the approximation error in terms of the best approximation error in two ways.

**Lemma 4.1.** *Let Assumption* (A3') *be satisfied. Then*

$$\mathcal{E}_{appr} \leq C_2 \mathcal{E}_{best} \, .$$

**Lemma 4.2.** *Let Assumption* (A3") *be satisfied. Then*

$$\mathcal{E}_{appr} \leq \frac{\Gamma}{2} \mathcal{E}_{best}^2 \, .$$

**Remark 4.3.** *It is worth noting the distinction between $\mathcal{E}_{appr}$ and $\mathcal{E}_{best}$. Both measure the minimal distance of elements $\Phi \in \mathcal{M}$ in the model class to the global minimizer $\Phi^*$. But $\mathcal{E}_{appr}$ measures this distance in terms of the cost function, while $\mathcal{E}_{best}$ measures it in the $\mathcal{V}$-norm.*

**Remark 4.4.** *The assumptions for both Lemmas are satisfied, for example, by linear and quadratic cost functionals, which are quite common in numerical applications.*

*Proof of Lemma 4.1.* Follows immediately. □

*Proof of Lemma 4.2.* We can express the functional by a first order Taylor expansion of the form

$$\mathcal{J}(\Phi) = \mathcal{J}(\Phi^*) + \mathrm{D}\mathcal{J}(\Phi^*)(\Phi - \Phi^*) + \frac{1}{2}\mathrm{D}^2\mathcal{J}(\xi)(\Phi - \Phi^*)(\Phi - \Phi^*) \, ,$$

for some $\xi \in \mathcal{V}$. Consequently, using $\Gamma = \sup_{\xi \in \mathcal{V}} \|\mathrm{D}^2\mathcal{J}(\xi)\|_{\mathcal{L}(\mathcal{V}, \mathcal{V}^*)} < \infty$ yields

$$|\mathcal{J}(\Phi) - \mathcal{J}(\Phi^*)| \leq \frac{\Gamma}{2} \|\Phi - \Phi^*\|_{\mathcal{V}}^2 \, .$$

This implies the claim. □

**Remark 4.5.** *Bounding $\mathcal{E}_{best}$ is an issue of approximation theory and for many model classes sharp bounds are known. However, approximation results for the popular deep neural networks, which could also be employed here, are scarce. This is the subject of ongoing research, see e.g. [35, 36].*

4.2. **Generalization error.** Motivated by the arguments and conclusions of Cucker and Smale [29] and Macdonald [37], the following analysis of the generalization error $\mathcal{E}_{\mathrm{gen}}$ relies on the concept of covering numbers, representing in a way the degree of compactness of the considered space. In our applications this seems to be a more natural notion than e.g. the classical VC dimension. For the sake of completeness, we provide full details of the derivations.

**Definition 4.6** (covering number). *The covering number $\nu(\mathcal{M}, \varepsilon)$ of a subset $\mathcal{M} \subseteq \mathcal{V}$ is defined as the minimal number of $\|\bullet\|_\mathcal{V}$-open balls of radius $\varepsilon$ needed to cover $\mathcal{M}$.*

**Example 4.7.** *Let $\mathcal{M}$ be a bounded subset of a finite dimensional linear space (e.g. a finite element space). The covering number of $\mathcal{M}$ can be estimated by*
$$\nu(\mathcal{M}, \varepsilon) \lesssim \mathrm{vol}(\mathcal{B}_1(0)) \left(\frac{R}{\varepsilon}\right)^{\dim(\mathcal{M})}.$$
*where $\mathcal{B}_1(0)$ denotes the unit ball and $R = \sup_{v \in \mathcal{M}} \|v\|$ is the radius of the domain.*

**Example 4.8.** *Let $\mathcal{M}$ be a Reproducing Kernel Hilber Space (RKHS). As sketched above, the minimization can be performed in a finite dimensional subspaces and the corresponding covering numbers are estimated in [30].*

**Example 4.9.** *Let $\mathcal{M}$ be a set of tensors in a given tensor format. Since $\mathcal{M}$ can be embedded into a finite dimensional linear space we obtain the bound from Example 4.7 as a crude upper bound. First sharper estimates for covering numbers for hierarchical tensors are provided by [38].*

**Example 4.10.** *For the currently extremely popular neural networks, estimates of the VC dimension (a quantity related to the covering number) are provided by [27] and first covering number estimates are derived in [36].*

**Remark 4.11.** *The concept of a covering number can be generalized to open balls with respect to a dissimilarity measure instead of the $\|\bullet\|_\mathcal{V}$ norm. This allows to apply the results of this section to more general loss functions $\ell$.*

The main convergence result for the generalization error is given in the following theorem.

**Theorem 4.12.** *For all $\varepsilon > 0$,*

$$\mathbb{P}[\mathcal{E}_{gen} > \varepsilon] \leq 2\nu(\mathcal{M}, \frac{\varepsilon}{8C_2})\delta(\tfrac{\varepsilon}{4}, N) \ ,$$

*where $\delta(\varepsilon, N)$ is an upper bound for $\mathbb{P}[|\mathcal{J}(\Phi) - \mathcal{J}_N(\Phi)| > \varepsilon]$.*

**Remark 4.13.** *Depending on $\nu(\mathcal{M}, \varepsilon)$ and the bound $\delta(\varepsilon, N)$, Theorem 4.12 may provide a very pessimistic upper bound with a relatively large pre-asymptotic range and a "phase shift" that separates an area of almost absolute certainty from an area where failure is almost guaranteed. This happens for example in high-dimensional linear spaces and when using the Hoeffding bound for $\mathbb{P}[|\mathcal{J}(\Phi) - \mathcal{J}_N(\Phi)| > \varepsilon]$. Nevertheless, the bound still provides a proof of convergence and illustrates the relation $N \in \mathcal{O}(\varepsilon^{-2}\ln(\varepsilon))$. This is slightly worse than the classical Monte Carlo bound of $N \in \mathcal{O}(\varepsilon^{-2})$ but is justified by the fact that we are not just evaluating integrals but are indeed optimizing a function $\Phi$ with respect to this integration. Moreover, we emphasize that the current analysis can not explain the significantly faster convergence we see in practical experiments as illustrated in Section 6. Better bounds may indeed be provided by exploiting further properties (e.g. sparsity) of the problem. This is illustrated in [12] for the linear case.*

We depict several lemmas in preparation of the proof of Theorem 4.12.

**Lemma 4.14.** *It holds,*

$$\mathcal{E}_{gen} \leq 2 \sup_{\Phi \in \mathcal{M}} |\mathcal{J}(\Phi) - \mathcal{J}_N(\Phi)|.$$

*Proof.* Recall that $\Phi^{*(\mathcal{M})}$ denotes a minimizer of $\mathcal{J}$ in the model class $\mathcal{M}$. We immediately derive

$$\begin{aligned}
\mathcal{E}_{\text{gen}} &= \mathcal{J}(\Phi^{*(\mathcal{M},N)}) - \mathcal{J}(\Phi^{*(\mathcal{M})}) \\
&= \mathcal{J}(\Phi^{*(\mathcal{M},N)}) - \mathcal{J}_N(\Phi^{*(\mathcal{M},N)}) + \mathcal{J}_N(\Phi^{*(\mathcal{M},N)}) - \mathcal{J}(\Phi^{*(\mathcal{M})}) \\
&\leq |\mathcal{J}(\Phi^{*(\mathcal{M},N)}) - \mathcal{J}_N(\Phi^{*(\mathcal{M},N)})| + \mathcal{J}_N(\Phi^{*(\mathcal{M})}) - \mathcal{J}(\Phi^{*(\mathcal{M})}) \\
&\leq 2 \sup_{\Phi \in \mathcal{M}} |\mathcal{J}(\Phi) - \mathcal{J}_N(\Phi)|.
\end{aligned}$$

□

**Lemma 4.15.** *Let $\varepsilon > 0$, $\nu := \nu(\mathcal{M}, \frac{\varepsilon}{8C_2})$ and $\{\Phi_j\}_{j \in [\nu]}$ the centers of the corresponding covering. Then it almost surely holds*

$$\sup_{\Phi \in \mathcal{M}} |\mathcal{J}(\Phi) - \mathcal{J}_N(\Phi)| < \frac{\varepsilon}{4} + \max_{1 \leq j \leq \nu} |\mathcal{J}(\Phi_j) - \mathcal{J}_N(\Phi_j)|.$$

*Proof.* Let $\Phi \in \mathcal{M}$ be given. By definition of the $\{\Phi_j\}_{j \in [\nu]}$ there exists some $\Phi_j$ with $\|\Phi - \Phi_j\|_{\mathcal{V}} < \frac{\varepsilon}{8C_2}$. Assumption (A2) implies $|\mathcal{J}(\Phi) - \mathcal{J}(\Phi_j)| < \frac{\varepsilon}{8}$ and almost surely $|\mathcal{J}_N(\Phi) - \mathcal{J}_N(\Phi_j)| < \frac{\varepsilon}{8}$. Hence,

$$
\begin{aligned}
&|\mathcal{J}(\Phi) - \mathcal{J}_N(\Phi)| \\
&\leq |\mathcal{J}(\Phi) - \mathcal{J}_N(\Phi) - (\mathcal{J}(\Phi_j) - \mathcal{J}_N(\Phi_j))| + |\mathcal{J}(\Phi_j) - \mathcal{J}_N(\Phi_j)| \\
&\leq |\mathcal{J}(\Phi) - \mathcal{J}(\Phi_j)| + |\mathcal{J}_N(\Phi) - \mathcal{J}_N(\Phi_j)| + |\mathcal{J}(\Phi_j) - \mathcal{J}_N(\Phi_j)| \\
&< \tfrac{\varepsilon}{4} + |\mathcal{J}(\Phi_j) - \mathcal{J}_N(\Phi_j)| \qquad \text{almost surely.} \\
&< \tfrac{\varepsilon}{4} + \max_{1 \leq j \leq \nu} |\mathcal{J}(\Phi_j) - \mathcal{J}_N(\Phi_j)| \qquad \text{almost surely.} \qquad \square
\end{aligned}
$$

**Lemma 4.16.** *Let $\varepsilon$, $\nu$ and $\{\Phi_j\}_{j \in [\nu]}$ be as in Lemma 4.15. Then,*

$$\mathbb{P}[\mathcal{E}_{gen} > \varepsilon] \leq \nu \max_{1 \leq j \leq \nu} \mathbb{P}[|\mathcal{J}(\Phi_j) - \mathcal{J}_N(\Phi_j)| > \tfrac{\varepsilon}{4}].$$

*Proof.* With the preceding lemmas, we deduce

$$
\begin{aligned}
\mathbb{P}[\mathcal{E}_{\text{gen}} > \varepsilon] &\leq \mathbb{P}[\sup_{\Phi \in \mathcal{M}} |\mathcal{J}(\Phi) - \mathcal{J}_N(\Phi)| > \tfrac{\varepsilon}{2}] &&\text{(Lemma 4.14)} \\
&\leq \mathbb{P}[\max_{1 \leq j \leq \nu} |\mathcal{J}(\Phi_j) - \mathcal{J}_N(\Phi_j)| > \tfrac{\varepsilon}{4}] &&\text{(Lemma 4.15)} \\
&\leq \sum_{1 \leq j \leq \nu} \mathbb{P}[|\mathcal{J}(\Phi_j) - \mathcal{J}_N(\Phi_j)| > \tfrac{\varepsilon}{4}] &&\text{(union bound)} \\
&\leq \nu \max_{1 \leq j \leq \nu} \mathbb{P}[|\mathcal{J}(\Phi_j) - \mathcal{J}_N(\Phi_j)| > \tfrac{\varepsilon}{4}]. \quad \square
\end{aligned}
$$

Using Lemma 4.16, the supremum can be factored out of the probability expression and then be bounded by means of classical concentration of measure arguments. This however comes at the price of the factor $\nu$. Two well-known upper bounds for the probability $\mathbb{P}[|\mathcal{J}(\Phi_j) - \mathcal{J}_N(\Phi_j)| > \tfrac{\varepsilon}{4}]$ are given by the Hoeffding and the Bernstein inequalities.

**Lemma 4.17** (Hoeffding 1963). *Let $\{X_i\}_{i=1,\ldots,N}$ be a sequence of i.i.d. bounded random variables $|X_i| \leq M$ and define $\overline{X} := \frac{1}{N} \sum_{i=1}^{N} X_i$. It then holds that*

$$\mathbb{P}\Big[|\mathbb{E}[\overline{X}] - \overline{X}| \geq \varepsilon\Big] \leq 2 \exp\left(-\frac{2\varepsilon^2 N}{M^2}\right).$$

**Lemma 4.18** (Bernstein 1927). *Let $\{X_i\}_{i=1,\ldots,N}$ be a sequence of i.i.d. bounded random variables $|X_i| \leq M$ with bounded variance $\mathrm{Var}(X_i) \leq \sigma^2$ and define $\overline{X} := \frac{1}{N} \sum_{i=1}^{N} X_i$. It then holds that*

$$\mathbb{P}\Big[|\mathbb{E}[\overline{X}] - \overline{X}| \geq \varepsilon\Big] \leq 2 \exp\left(-\frac{\tfrac{1}{2}\varepsilon^2 N}{\sigma^2 + \tfrac{1}{3} M \varepsilon}\right).$$

**Corollary 4.19.** *If Assumptions (A1) holds, then Lemma 4.17 leads to the estimate*

$$\delta(\varepsilon, N) \leq 2 e^{-2\varepsilon^2 N / C_1^2} .$$

*If the variance $\sigma^2$ of $\ell(\Phi, \bullet)$ can be assumed to be negligible, Lemma 4.18 yields the even tighter bound*

$$\delta(\varepsilon, N) \leq 2e^{-3\varepsilon N/(4C_1^2)} \ .$$

*Proof.* Recall that the samples $x^i$ are i.i.d. and therefore the $\ell(\Phi, x^i)$ are i.i.d. for $i = 1, \ldots, N$. Moreover, (A1) ensures that $|\ell(\Phi, x^i)| \leq C_1$ holds almost surely for all $i = 1, \ldots, N$. Hence the assumptions for Lemma 4.17 are satisfied. Finally, (A1) also ensures that $\mathrm{Var}(\ell(\Phi, x^i)) \leq C_1^2$ holds almost surely for all $i = 1, \ldots, N$. Therefore, Lemma 4.18 is applicable. □

Theorem 4.12 now is a mere corollary of the preceding lemmas.

*Proof of Theorem 4.12.* Follows from Lemma 4.16 and Corollary 4.19. □

4.3. **Norm Error.** We aim at using the proposed framework in the numerical analysis of parametric PDEs as outlined in Section 5. Hence, we are primarily interested in the error of the parameters $\mathcal{E}_{\mathrm{norm}} = \|\Phi^* - \Phi^{*(\mathcal{M},N)}\|_{\mathcal{V}}$. This section is devoted to the derivation of respective error bounds.

**Lemma 4.20.** *Let Assumption* (A4) *on local strong convexity be satisfied and assume that $\Phi^{*(\mathcal{M},N)}$ lies in the strongly convex neighbourhood $\mathcal{U}$ of $\Phi^*$. Then*

$$\mathcal{E}_{norm}^2 \leq \frac{2}{\gamma}\mathcal{E}_{cost} \ .$$

*Proof.* By the first order optimality condition for $\Phi^*$,

$$D\mathcal{J}(\Phi^*)(\Phi - \Phi^*) = 0 \quad \text{for all } \Phi \in \mathcal{V}.$$

Consequently, by Assumption (A4),

$$\begin{aligned}
\mathcal{E}_{\mathrm{norm}}^2 &= \|\Phi^* - \Phi^{*(\mathcal{M},N)}\|^2 \\
&\leq \frac{2}{\gamma}\Big(\mathcal{J}(\Phi^{*(\mathcal{M},N)}) - \mathcal{J}(\Phi^*) - D\mathcal{J}(\Phi^*)(\Phi^{*(\mathcal{M},N)} - \Phi^*)\Big) \\
&= \frac{2}{\gamma}|\mathcal{J}(\Phi^*) - \mathcal{J}(\Phi^{*(\mathcal{M},N)})| \ .
\end{aligned}$$

This proves the claim. □

**Lemma 4.21.** *Let Assumptions* (A3") *and* (A4) *be satisfied. Then*

$$\mathcal{E}_{norm}^2 \leq \frac{\Gamma}{\gamma}\mathcal{E}_{best}^2 + \frac{2}{\gamma}\mathcal{E}_{gen} \ .$$

An interpretation of this estimate is that $\mathcal{E}_{\text{gen}}$ in a way measures the quality of the chosen cost functional. If the functional is chosen poorly, the minimum of the empirical functional $\Phi^{*(\mathcal{M},N)}$ may deviate strongly from the minimum $\Phi^{*(\mathcal{M})}$ and may only converge with a slow rate.

*Proof of Lemma 4.21.* By Lemma 4.20, the splitting of $\mathcal{E}_{\text{cost}}$ and Lemma 4.2,

$$\begin{aligned}\mathcal{E}_{\text{norm}}^2 &\leq \frac{2}{\gamma}\mathcal{E}_{\text{cost}} \\ &\leq \frac{2}{\gamma}(\mathcal{E}_{\text{appr}} + \mathcal{E}_{\text{gen}}) \\ &\leq \frac{2}{\gamma}\left(\frac{\Gamma}{2}\mathcal{E}_{\text{best}}^2 + \mathcal{E}_{\text{gen}}\right).\end{aligned}$$ □

**Corollary 4.22.** *Let $a > 0$ and define*

$$p(a, N) := \mathbb{P}[\mathcal{E}_{gen} \geq \tfrac{1}{2}a\Gamma\mathcal{E}_{best}^2] .$$

*Then, with probability $1 - p(a, N)$,*

$$\mathcal{E}_{best}^2 \leq \mathcal{E}_{norm}^2 \leq (1+a)\frac{\Gamma}{\gamma}\mathcal{E}_{best}^2 . \tag{1}$$

Recall that the norm error $\mathcal{E}_{\text{norm}}$ is the distance of the empirical solution $\Phi^{*(\mathcal{M},N)}$ to the exact solution $\Phi^*$ and that the best approximation error $\mathcal{E}_{\text{best}}$ is the smallest such distance that is obtainable in the given model class. Corollary 4.22 shows that in a "validity probability" the norm error is quasi-optimal in the sense that it is equivalent to the best approximation error. Using Lemma 4.12 we can bound $p(a, N)$ in Corollary 4.22 and show that the probability $1 - p(a, N)$ of equivalence (1) tends to one with an exponential rate in the number of samples $N$. This means that the convergence results hold with high probability, provided that $N$ is sufficiently large.

*Proof of Corollary 4.22.* Lemma 4.21 implies

$$\begin{aligned}\mathbb{P}[\mathcal{E}_{\text{norm}}^2 \geq \varepsilon] &\leq \mathbb{P}\left[\frac{\Gamma}{\gamma}\mathcal{E}_{\text{best}}^2 + \frac{2}{\gamma}\mathcal{E}_{\text{gen}} \geq \varepsilon\right] \\ &\leq \mathbb{P}\left[\mathcal{E}_{\text{gen}} \geq \frac{1}{2}(\gamma\varepsilon - \Gamma\mathcal{E}_{\text{best}}^2)\right] .\end{aligned}$$

Choosing $\varepsilon = \gamma^{-1}(1+a)\Gamma\mathcal{E}_{\text{best}}^2$ yields the assertion. □

## 5. Problems in UQ

In this section we discuss the intended application of the proposed method in the context of Uncertainty Quantification. More specifically,

consider the abstract problem
$$\mathcal{D}(u; y) = 0\,,$$

where $\mathcal{D}$ encodes a (possibly non-linear) PDE model in a physical domain $D \subset \mathbb{R}^d$, $d = 1, 2, 3$, depending on parameters $y \in \Gamma \subset \mathbb{R}^M$. The parameter vector $y$ may be finite dimensional ($M < \infty$) or infinite dimensional ($M = \infty$). We assume that the parametric solution $u(x, y)$ can be represented as $u \in L^2(\Gamma, \rho; \mathcal{X}) \simeq \mathcal{X} \otimes \mathcal{Y} =: \mathcal{V}$ with (typically) $\mathcal{X} \subseteq H^1(D)$ and $\mathcal{Y} \subseteq L^2(\Gamma, \rho)$ and $\rho$ some probability measure on $\Gamma$.

Further details of these problems can e.g. be found in [6–8, 10, 20]. Moreover, adaptive Galerkin discretizations are for instance considered in [39–43], which is one of the UQ standard methods we have in mind.

5.1. **Parametric PDEs.** As a common benchmark example, we introduce a linear elliptic PDE with homogeneous Dirichlet boundary data where the solution $u(x, y)$ solves
$$-\nabla \cdot (a(x, y) \nabla u(x, y)) = f(x), \qquad u|_{\partial D} = 0. \tag{2}$$

This equation models the stationary density $u(x, y)$ of a substance diffusing through a medium with permeability $a(x, y) > 0$. The parametric coefficient is assumed to be given either by an affine representation of the form
$$a(x, y) = a_0(x) + \sum_{m=1}^{M} \sigma_m a_m(x) y_m \quad \text{with} \quad y \in \Gamma = [-1, 1]^M, \tag{3}$$

or by the numerically more involved representation
$$a(x, y) = \exp\left(a_0(x) + \sum_{m=1}^{M} \sigma_m a_m(x) y_m\right) \quad \text{with} \quad y \in \Gamma = \mathbb{R}^M\,. \tag{4}$$

In both cases $\{a_0, a_1, \ldots, a_M\}$ is considered to be an orthonormal basis in $L^2(D)$ and the coefficients $\sigma_m > 0$ are assumed to be positive. The parameter $y \in \Gamma$ can be associated with a random variable $y \sim \rho$ that determines uncertainty in the porosity of the medium. For the affine case (3) $\rho$ is chosen to be a uniform distribution $\rho = \mathcal{U}(\Gamma)$ while for the log-normal case (4) it is chosen to be standard Gaussian $\rho = \mathcal{N}(0, I)$. We assume the problem to be elliptic with high probability (uniform ellipticity assumption), i.e., for a small $\varepsilon > 0$, there exist constants $\underline{a}, \overline{a} > 0$ such that
$$P\left[0 < \underline{a} \leq a(x, y) \leq \overline{a} < \infty\right] > 1 - \varepsilon \qquad \text{a.s. for } (x, y) \in D \times \Gamma.$$

We point out that for most random fields, choosing a finite parameter dimension $M$ is a required simplifying restriction for the actual computation to become feasible. Adaptive a priori or a posteriori methods are available to sensibly control this parameter. For the Stochastic Galerkin method, the affine setting is e.g. examined in [39, 41, 44] and first results with log-normal coefficients are shown in [42] for reliable residual-based error estimators. Other numerical methods such as Stochastic Collocation [45, 46] rely on a piori estimators or heuristic hierarchical indicators.

We employ the variational formulation of problem (2) with respect to a test function $v \in \mathcal{V}$, cf. [39, 41]. It reads

$$A(u,v) := \int_\Gamma \int_D a(x,y) \nabla u(x,y) \cdot \nabla v(x,y) \rho(y) \, \mathrm{d}x \, \mathrm{d}y$$
$$= \int_\Gamma \int_D f(x,y) v(x,y) \rho(y) \, \mathrm{d}x \, \mathrm{d}y \ .$$

Since $\mathcal{V}$ is a Hilbert space, this variational problem can be treated directly by means of the presented *Variational Monte Carlo* approach by defining the bilinear form $A(u,v)$, which induces $\|v\|_A^2 := A(v,v)$. We then consider the equivalent optimization problem

$$u = \underset{\Phi \in \mathcal{V}}{\arg\min} \, \frac{1}{2} \|A\Phi - f\|_\mathcal{V}^2$$
$$= \underset{\Phi \in \mathcal{V}}{\arg\min} \, \frac{1}{2} \|A(\Phi - u)\|_\mathcal{V}^2 \ .$$

This means that the solution $u \in \mathcal{V}$ can be obtained by minimizing

$$\mathcal{J}(\Phi) := \frac{1}{2} \int_\Gamma \|\Phi(y) - u(y)\|_{B(y)}^2 \rho(y) \, \mathrm{d}y \ ,$$

with the bilinear form $B(y) \colon \mathcal{X} \times \mathcal{X} \to \mathbb{R}$ defined by

$$B(y)(u,v) := \int_D a(x,y) \nabla u(x) \nabla v(x) \, \mathrm{d}x \ ,$$

and $\|v\|_{B(y)}^2 := B(y)(v,v)$.

To find an appropriate finite dimensional model class, recall that $\mathcal{V} = H^1(D) \otimes L^2(\Gamma, \rho)$. We first choose finite dimensional subspaces

$$\mathcal{X}_{\text{FEM}} \subseteq \mathcal{X} \quad \text{and} \quad \mathcal{Y}_\mathcal{I} \subseteq \mathcal{Y},$$

and consider the *discrete solution subspace*

$$\mathcal{V}_\mathcal{I} = \mathcal{X}_{\text{FEM}} \otimes \mathcal{Y}_\mathcal{I} \subseteq \mathcal{V} \ ,$$

generated explicitly by the $S$-dimensional conforming FE space $\mathcal{X}_{\text{FEM}} = \text{span}\{\varphi_j\}_{j=1,\ldots,S}$ with respect to a regular triangulation of the domain $D$. Since $\mathcal{Y} = L^2(\Gamma, \rho)$ exhibits a (countable) product structure, we can choose $\mathcal{Y}_\mathcal{I}$ as the vector space generated by a tensor product basis of multivariate polynomials with $\alpha \in \mathcal{I}$ given by

$$P_\alpha(y) = \bigotimes_{m=1}^{M} P_{\alpha_m}(y_m).$$

Here, the $P_{\alpha_m}$ are orthogonal w.r.t. the marginal distribution $\rho_m$ and $\mathcal{I}$ is a finite subset of $\prod_{m=1}^{M}[q_m]$ for given $q_m \in \mathbb{N}$. Hence, every $\Phi \in \mathcal{V}_\mathcal{I}$ can be written as

$$\Phi(x,y) = \Phi_W(x,y) = \sum_{j=1}^{S} \sum_{\alpha \in \mathcal{I}} W(j,\alpha) \varphi_j(x) P_\alpha(y),$$

with coefficient tensor $W \in \mathbb{R}^{S \times q_1 \times \cdots \times q_M}$. As described in [1, 41, 47–49], we can represent $W$ by

$$W(j,\alpha) = U_0(j) U_1(\alpha_1) \cdots U_M(\alpha_M),$$

with

$$\begin{aligned} U_0(j) &\in \mathbb{R}^{r_0}, & j &= 1, \ldots, S, \\ U_i(\alpha_i) &\in \mathbb{R}^{r_{i-1} \times r_i}, & \alpha_i &= 1, \ldots, q_i, \\ U_M(\alpha_M) &\in \mathbb{R}^{r_{M-1}}, & \alpha_M &= 1, \ldots, q_M, \end{aligned}$$

for some *ranks* $r_0, \ldots, r_{M-1} \in \mathbb{N}$. The set of all coefficient tensors $W$ with prescribed ranks $\mathbf{r} = (r_0, \ldots, r_{M-1})$ is denoted by $\mathcal{M}_{\leq r}$. We now define the model class by

$$\mathcal{M} = \{\Phi_W \in \mathcal{V}_\mathcal{I} : W \in \mathcal{M}_{\leq r}\}.$$

Denote by $u_h(y)$ the FE solution of $B(y)u_h(y) = f(y)$ in $\mathcal{X}_{\mathrm{FEM}}$. Then by Galerkin orthogonality it holds that for all $\Phi \in \mathcal{M}$ that

$$\begin{aligned} \mathcal{J}(\Phi) &= \int_\Gamma \|\Phi(y) - u(y)\|_{B(y)}^2 \rho(y) \, \mathrm{d}y \\ &= \int_\Gamma \|\Phi(y) - u_h(y)\|_{B(y)}^2 \rho(y) \, \mathrm{d}y. \end{aligned}$$

Due to the uniform ellipticity assumption, the norms $\|\bullet\|_{B(y)}$ and $\|\bullet\|_{H_0^1(D)}$ are equivalent with high probability, i.e.,

$$\underline{a} \|\bullet\|_{H_0^1(D)} \leq \|\bullet\|_{B(y)} \leq \overline{a} \|\bullet\|_{H_0^1(D)}.$$

We thus can introduce the functional

$$\widetilde{\mathcal{J}}(\Phi) = \int_\Gamma \overline{a}^2 \|\Phi(y) - u_h(y)\|_{H_0^1(D)}^2 \rho(y) \, \mathrm{d}y,$$

which is numerically much easier to handle than $\mathcal{J}(\Phi)$. The corresponding minimum $\widetilde{\Phi} \in \mathcal{M}$ is quasi-optimal by the equivalence of the norms,

namely

$$\int_\Gamma \|\widetilde{\Phi}(y) - u(y)\|^2_{B(y)} \rho(y) \, \mathrm{d}y \leq \int_\Gamma \overline{a}^2 \|\widetilde{\Phi}(y) - u(y)\|^2_{H^1_0(D)} \rho(y) \, \mathrm{d}y$$

$$\leq \int_\Gamma \overline{a}^2 \|\Phi^{*(\mathcal{M})}(y) - u(y)\|^2_{H^1_0(D)} \rho(y) \, \mathrm{d}y$$

$$\leq \int_\Gamma \frac{\overline{a}^2}{\underline{a}^2} \|\Phi^{*(\mathcal{M})}(y) - u(y)\|^2_{B(y)} \rho(y) \, \mathrm{d}y \, .$$

Consequently, it follows with high probability that

$$\|\widetilde{\Phi} - u\|_\mathcal{V} \leq \frac{\overline{a}}{\underline{a}} \|\Phi^{*(\mathcal{M})} - u\|_\mathcal{V} \quad \text{and} \quad \|\widetilde{\Phi} - u\|_A \leq \frac{\overline{a}}{\underline{a}} \|\Phi^{*(\mathcal{M})} - u\|_A \, .$$

Finally, in accordance with the proposed framework, this is formulated as empirical functional

$$\widetilde{\mathcal{J}}_N(\Phi) := \frac{1}{N} \sum_{i=1}^N \|\Phi(y^i) - u(y^i)\|^2_{H^1_0(D)} \, , \tag{5}$$

with $y^i$ sampled according to the density $\rho$. The approximations $u_h(y^i)$ can be obtained numerically by standard FE methods.

In the present setting, $a(\bullet, y)$ and $y$ are in a one-to-one relation. Thus, the central notion of this scheme is to **learn the solution operator** $y \mapsto a(\bullet, y) \mapsto u(\bullet, y)$ from generated data $u(\bullet, y^i)$. The computation of the data $u(\bullet, y^i)$ is completely non-intrusive, meaning that standard FE solvers can be employed to compute $u(\bullet, y^i)$ and perform the optimization of the mean squared errors in a standard tensor recovery algorithm as described in [50].

**Remark 5.1.** *For numerical reasons it is beneficial to represent the samples in a basis $\Psi$ that is orthogonal w.r.t. the $H^1_0(D)$ scalar product. If $\Phi$ denotes the vector of standard FE basis functions, every function $u \in \mathcal{X}_{FEM}$ can be represented by*

$$u(x) = \overline{u}^T \Phi(x) = \underline{u}^T \Psi(x) \, ,$$

*where $\overline{u}$ and $\underline{u}$ denote the coefficient vectors. Using the stiffness matrix $S$, the $H^1_0(D)$ scalar product can then be computed via*

$$(u, v)_{H^1_0(D)} = \overline{u}^T S \overline{v} = \underline{u}^T \underline{v}.$$

*This makes obvious that any decomposition $S = X^T X$ provides a basis transform $X$ from $\Phi$ to a basis $\Psi = X^{-T} \Phi$ satisfying the orthogonality condition. In fact, this allows for the computation of the $H^1_0(D)$ norm of functions in the FE space more efficiently and provides numerical stability to the minimization problem* (5). *Note also that one can use a standard tensor reconstruction algorithm as in [50] since*

$\|\underline{u}\|_{H_0^1(D)} = \|\underline{u}\|_{\ell^2}$. *However, the resulting tensor represents the solution w.r.t. the orthogonal basis $\Psi$ and another basis transformation after optimization is thus required. Moreover, the price for the described efficiency and numerical stability is a costly and numerically unstable Cholesky factorization of the ill-conditioned stiffness matrix $S$.*

5.2. **Backward Kolmogorov and Fokker-Planck equations.** As another viable application of our Variational Monte Carlo method, we derive a suitable formulation of the Backward Kolmogorov and Fokker-Planck equations, which are high-dimensional elliptic PDEs derived from stochastic processes. Let us consider a Langevin dynamic driven by a gradient field with a *confining potential* $V : \mathbb{R}^d \to \mathbb{R}$, $V \in C^2(\mathbb{R}^d)$, such that $V(x) \to \infty$ for $|x| \to \infty$ and $x \mapsto e^{-V(x)} \in L^1(\mathbb{R}^d)$. We then define the following stochastic differential equations subject to a stochastic process $\{X_t\}_{t\geq 0}$ with a $d$-dimensional random variable $X_t$,

$$\mathrm{d}X_t = -\mathrm{grad}\, V(X_t)\,\mathrm{d}t + \sigma\,\mathrm{d}W_t\ ,$$

where $\{W_t\}_{t\geq 0}$ is a standard $d$-dimensional Wiener process [51]. The probability density $p$ of finding a particle at position $x$ and time $t$ is governed by the *Fokker-Planck (FP) equation*

$$\frac{\partial}{\partial t}p = L^\dagger p := \mathrm{div}((\mathrm{grad}\, V)p) + \frac{\sigma^2}{2}\Delta p\ .$$

The formal adjoint of the Fokker-Planck operator leads to the Backward Kolmogorov equation

$$\frac{\partial}{\partial t}u = Lu := -\mathrm{grad}\, V \cdot \mathrm{grad}\, u + \frac{\sigma^2}{2}\Delta u\ .$$

The presented approach can be used to solve an implicit Euler step for the Backward Kolmogorov equation of the form

$$(I - hL)u_{k+1} = u_k\ , \tag{6}$$

with $u_k(x) := u(t_k, x)$ and $t_k := hk$.

**Remark 5.2.** *The solution of the corresponding Fokker-Planck equation is given by $p(t,x) = \rho(x)u(t,x)$ and can thus be computed indirectly by the Variational Monte Carlo method.*

Under suitable conditions [52], to solve (6) we consider the equation in a weighted $L^2$ space. The weight is chosen as the equilibrium density

$$\rho(x) := \frac{1}{Z}e^{-\frac{2}{\sigma^2}V(x)}\ ,\ x \in \mathbb{R}^d\ ,$$

with the normalization constant $Z = \int_{\mathbb{R}^d} e^{-\frac{2}{\sigma^2}V(x)}\,\mathrm{d}x$ which satisfies the stationary equation $L^\dagger \rho = 0$. The backward Kolmogorov operator $L : H^1(\mathbb{R}^d, \rho) \to H^1(\mathbb{R}^d, \rho)$ is considered as acting on the weighted Sobolev space

$$H^1(\mathbb{R}^d, \rho) := \{u \in L^2(\mathbb{R}^d, \rho) \; : \; \nabla u \in L^2(\mathbb{R}^d, \rho; \mathbb{R}^d)\} \;,$$

with the inner product

$$a(u, v) := \langle \nabla u, \nabla v \rangle_{L^2(\mathbb{R}^d, \rho; \mathbb{R}^d)} + \langle u, v \rangle_{L^2(\mathbb{R}^d, \rho)} \;.$$

It is straightforward to verify that in this weighted space the backward Kolmogorov operator $L$ becomes symmetric [52] and that solving equation (6) is equivalent to the minimization of the cost functional

$$\mathcal{J}(v) := \frac{1}{2}(\langle v, v \rangle_{L^2(\mathbb{R}^d, \rho)} + h\frac{\sigma^2}{2}\langle \nabla v, \nabla v \rangle_{L^2(\mathbb{R}^d, \rho; \mathbb{R}^d)}) - \langle u_k, v \rangle_{L^2(\mathbb{R}^d, \rho)} \;. \quad (7)$$

**Theorem 5.3.** *For all $u, v \in H_0^1(\mathbb{R}^d, \rho)$ it holds that*

$$\langle u, Lv \rangle_{L^2(\mathbb{R}^d, \rho)} = -\frac{\sigma^2}{2} \langle \nabla u, \nabla v \rangle_{L^2(\mathbb{R}^d, \rho; \mathbb{R}^d)} \;.$$

*As a consequence, $L$ is symmetric w.r.t. the prescribed weighted $L^2$ inner product.*

*Proof.* In the following we denote by $\langle \bullet, \bullet \rangle$ the $L^2$ inner product and by $\langle \bullet, \bullet \rangle_\rho$ the weighted $L^2$ inner product w.r.t. the weight $\rho$. Integration by parts and the product rule yield

$$\begin{aligned}\langle u, Lv \rangle_\rho &= \frac{\sigma^2}{2}\langle u, \Delta v \rangle_\rho - \langle u, (\nabla V) \cdot (\nabla v) \rangle_\rho \\ &= -\frac{\sigma^2}{2}\langle \nabla(u\rho), \nabla v \rangle - \langle u, (\nabla V) \cdot (\nabla v) \rangle_\rho \\ &= -\frac{\sigma^2}{2}\langle \nabla u, \nabla v \rangle_\rho - \frac{\sigma^2}{2}\langle u\nabla\rho, \nabla v \rangle - \langle u, (\nabla V) \cdot (\nabla v) \rangle_\rho \;.\end{aligned}$$

Now observe that by definition of $\rho$,

$$\nabla \rho = -\frac{2}{\sigma^2}(\nabla V)\rho \;.$$

Substitution concludes the proof. $\square$

Clearly, we can use the Variational Monte Carlo framework to minimize the cost functional (7). This provides a further striking example for the proposed methodology. The solution of the eigenvalue problem for the backward Kolmogorov operators including more details about the operators is considered in a forthcoming paper, which will also contain numerical verifications of the suggested approach. Another way to obtain a solution of the backward Kolmogorov equation by deep neural networks is considered in a recent publication [53].

6. Numerical Experiments

This section is concerned with illustrating the performance of the proposed Variational Monte Carlo approach for a set of benchmark problems with parametric PDEs. To assess the quality of the resulting approximation of the parametric solution, the error of the expectation and the variance with respect to a Quasi-Monte Carlo (QMC) reference solution from $10^6$ Sobol samples is depicted.

The used model classes and cost functionals are defined analogously to Section 5. In the resulting tensor representation, the high-dimensional integrals (quantities of interest such as the expectation) can be evaluated efficiently and exactly. The computation of the FE solution samples is carried out with a standard conforming $P_1$-Galerkin method on the domain $D = (0,1)^2$ using the open source `FEniCS` software package [54]. For the tensor reconstructions the Block-ASD optimization algorithm described in [50] is employed. A complete implementation is contained in the open source `C++` library `xerus` [55], which is interfaced in the open source `python` framework `ALEA` [56]. For the reconstruction, a maximal polynomial degree of 12 (affine) and 4 (log-normal) is allowed for each polynomial basis in the stochastic dimensions. We plot the error progression with respect to the number of reconstruction samples. The rank is chosen adaptively with an allowed maximum of 40, which for most problems does not impose a restriction. We prevent excessive rank increases by allowing each rank to increase only once every 10 iterations. For the experiments, we follow and extend our previous work [50] and use the same optimization algorithm for the tensor reconstruction.

6.1. **Parametric PDEs.** We consider parametric second order PDEs in the context of UQ as described in Section 5. The model for the coefficient $a(x,y)$ is given by

$$a(x,y) := a_0(x) + \sum_{m=1}^{M} \sigma_m a_m(x) y_m.$$

Here, $a_0$ determines the mean field, which is set to 1 if not stated otherwise, and $a_m(x) \propto m^{-2} \sin(\lfloor \frac{m+2}{2} \rfloor \pi x_1) \sin(\lceil \frac{m+2}{2} \rceil \pi x_2)$. This choice corresponds to the slow decay experiments e.g. in [39, 40] but with a larger scaling factor of the modes.

The examined examples are defined by the following PDEs with homogeneous Dirichlet boundary conditions and $a(x,y)$ expanded in $M = 20$ terms[1]. The relative errors for the expectation value and variance obtained by the reconstruction are compared to Monte Carlo and Quasi Monte Carlo simulations in Figures 1-4.

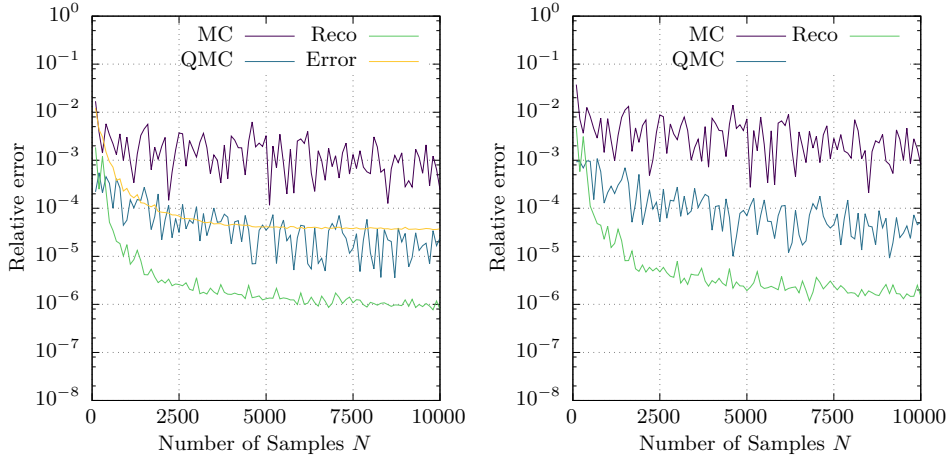

**Figure 1.** PDE setting (I): TT reconstruction error of the *expectation* ("Reco") compared to MC and QMC simulations (left). Additionally the average relative error for a random set of 1000 samples *not* used for the reconstruction is shown ("Error"). TT reconstruction error of the *variance* compared to MC and QMC simulations (right).

(I) **Diffusion** (affine)
$$-\nabla \cdot (a(x,y)\nabla u(x,y)) = 1 \quad \text{with} \quad y \sim \mathcal{U}([-1,1]^M)$$

(II) **Diffusion** (lognormal, $a_0 \equiv 0$)
$$-\nabla \cdot (\exp(a(x,y))\nabla u(x,y)) = 1 \quad \text{with} \quad y \sim \mathcal{N}(0, I_M)$$

(III) **Nonlinear Diffusion** (affine)
$$-\nabla \cdot \left(\left(\frac{a(x,y)}{10} + u(x,y)\right)^2 \nabla u(x,y)\right) = 1 \quad \text{with} \quad y \sim \mathcal{U}([-1,1]^M)$$

(IV) **Convection-Diffusion** (affine, SUPG stabilized FEM)
$$-\nabla \cdot (\kappa \nabla u(x,y)) + \beta \cdot \nabla u(x,y) = 1,$$
with $\kappa = 10^{-2}$ and $\beta = \begin{pmatrix} 1 - a(x,y) & 1 - |a(x,y)| \end{pmatrix}^\mathsf{T}$.

The examples demonstrate that the suggested approach can be sucessfully applied to linear and non-linear PDE problems alike. In Figures 1

---

[1]Note that the parameter vector $y$ is the image of a random variable and hence is associated to a probability distribution but (despite the inaccurate notation) not actually a random quantity.

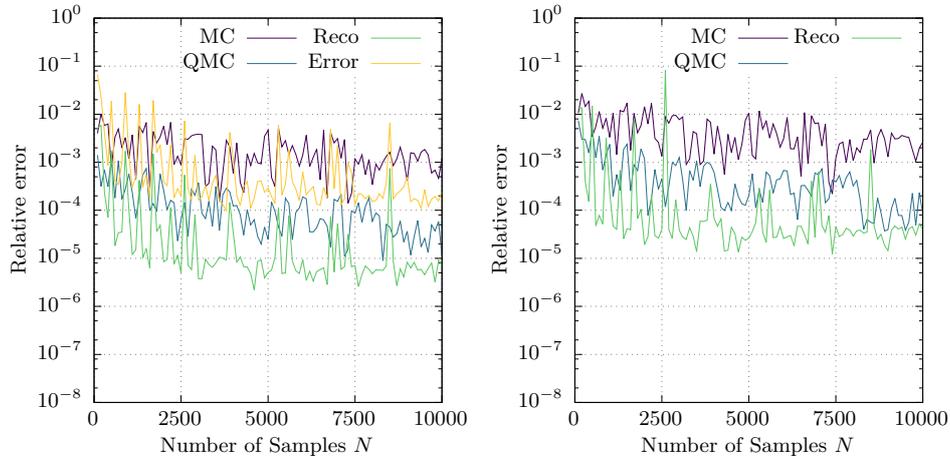

**Figure 2.** PDE setting (II): TT reconstruction error of the *expectation* ("Reco") compared to MC and QMC simulations (left). Additionally the average relative error for a random set of 1000 samples *not* used for the reconstruction is shown ("Error"). TT reconstruction error of the *variance* compared to MC and QMC simulations (right).

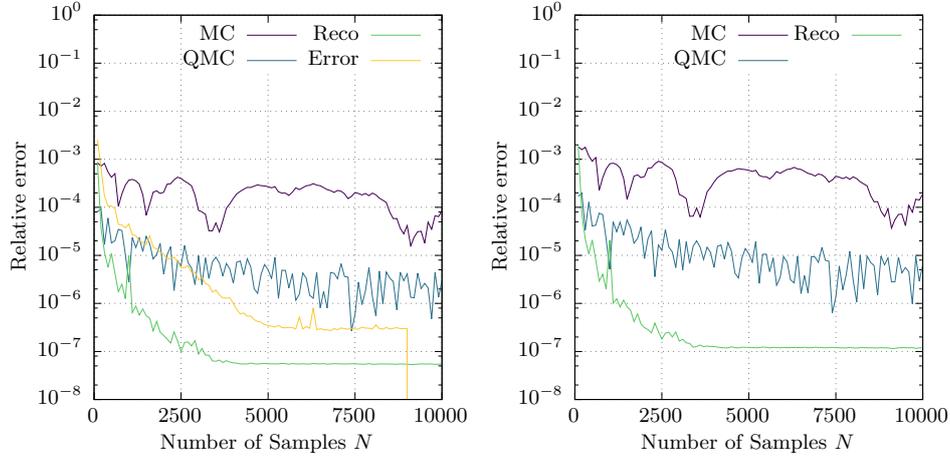

**Figure 3.** PDE setting (III): TT reconstruction error of the *expectation* ("Reco") compared to MC and QMC simulations (left). Additionally the average relative error for a random set of 1000 samples *not* used for the reconstruction is shown ("Error"). TT reconstruction error of the *variance* compared to MC and QMC simulations (right).

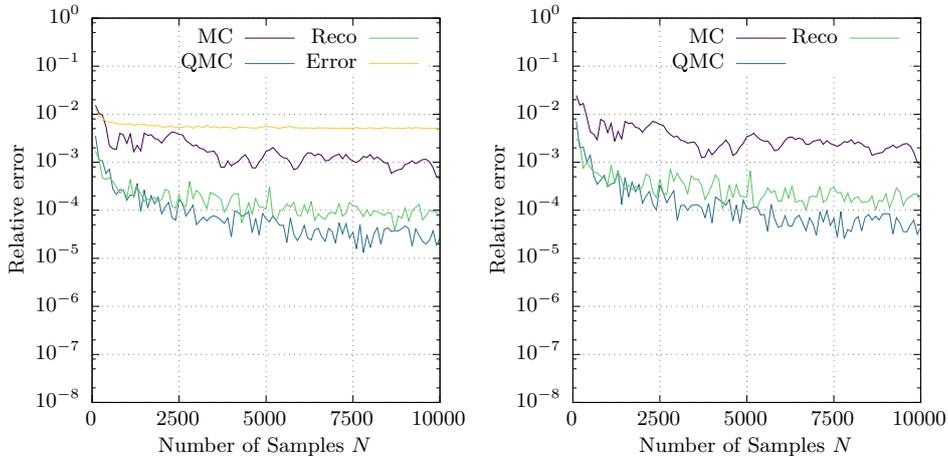

**Figure 4.** PDE setting (IV): TT reconstruction error of the *expectation* ("Reco") compared to MC and QMC simulations (left). Additionally the average relative error for a random set of 1000 samples *not* used for the reconstruction is shown ("Error"). TT reconstruction error of the *variance* compared to MC and QMC simulations (right).

to 4, the relative errors for the expectation and variance obtained by the tensor reconstruction are compared to Monte Carlo simulations. We emphasize again that the tensor reconstruction in fact represents the entire stochastic parametric solution $u(x,y) \in \mathcal{V}$ from which the first two moments are evaluation in order to compare the approximation quality with classical methods, which only allow for an evaluation of functionals. As a consequence, in principle arbitrary statistical quantities can be obtained globally and locally from the tensor representation of $u(x,y)$ (e.g. pointwise densities).

In all cases we observe that the presented approach provides significantly better results for the expectation than Monte Carlo sampling. It is striking that for the Darcy examples, this even holds when compared with the QMC simulations. The log-normal case (II) is considered as rather involved and the observed results are quite ecouraging. We note that the pointwise error ("Error"), which represents the approximation quality of the actual parametric solution (not just a functional), is also small.

For the complicated settings (III) and (IV) with either non-linear $x$-dependence or explicit non-linear $y$-dependence, the QMC results are slightly better for the first two moments when compared to the tensor reconstruction. Nevertheless, the relative errors are of the same order of magnitude. A good indication for the large variance of the solution manifold is the now larger pointwise "Error". However, it should be noted that this is still quite accurate for the types of problems considered here. It can be expected that a sensible increase of the tensor ranks, the polynomial degree and the number of samples used for the reconstruciton would lead to even better results.

6.2. **Cookie Problems.** In this section, in order to examine a different type of setup with inherently finite-dimensional noise, we consider two so-called "cookie problems". With these, circular inclusions of fixed or random size and with different random diffusion coefficients are prescribed.

(V) **Diffusion** (fixed radii) Let 9 subdomains of $D$ be given by discs $D_k$ ($k = 1, \ldots, 9$) with fixed radius $r = 1/8$ and centers $c = \begin{pmatrix} i/6 & j/6 \end{pmatrix}^\intercal$ for $i, j \in \{1, 3, 5\}$. The considered problem depends on $y = \begin{pmatrix} y_1 & \ldots & y_9 \end{pmatrix}^\intercal$ with $y_k \sim \mathcal{U}(-1, 1)$, has homogeneous Dirichlet boundary conditions and is given by

$$-\nabla \cdot (\kappa(x,y)\nabla u(x,y)) = 1,$$

where $\kappa|_{D\setminus \cup_{k=1,\ldots,4} D_k} = 1$ and $\kappa|_{D_k} = y_k$.

(VI) **Diffusion** (random radii) The setting is the same as before. However, the problem depends on additional parameters $y^* = \begin{pmatrix} y_{10}^* & \ldots & y_{18}^* \end{pmatrix}^\intercal$ with $y_j^* \sim \mathcal{U}(-1, 1)$ determining the radii $r_j = (1 + y_j^*/3)/3$ of the inclusions $j = 10, \ldots, 18$.

The relative errors for the expectation value and variance obtained by the reconstruction are compared to (Quasi) Monte Carlo simulations in Figures 5-6. As before, the tensor reconstruction yields significantly more accurate results (fixed radii) or at least errors comparable to the QMC simulations. This is especially noteworthy since setting (VI) clearly exhibits no tensor structure. When compared to the simpler PDE setting (V), all depicted results apparently indicate a much more difficult problem.

## 7. Assessment and Outlook

The analysis in this paper is based on known results from statistical learning, which we scrutinize from the perspective of UQ and high-dimensional PDEs. Central references for our exposition are [29, 30] with additional relations to the theory of Chervonenkis & Vapnik. We consider the presented treatise as an initial step in this direction with rigorous error estimates but with many remaining open questions.

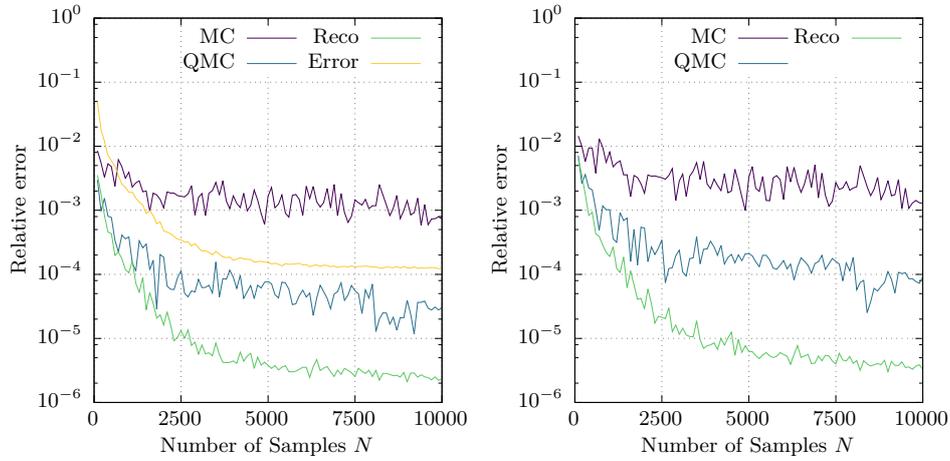

**Figure 5.** PDE setting (V): TT reconstruction error of the *expectation* ("Reco") compared to MC and QMC simulations (left). Additionally the average relative error for a random set of 1000 samples *not* used for the reconstruction is shown ("Error"). TT reconstruction error of the *variance* compared to MC and QMC simulations (right).

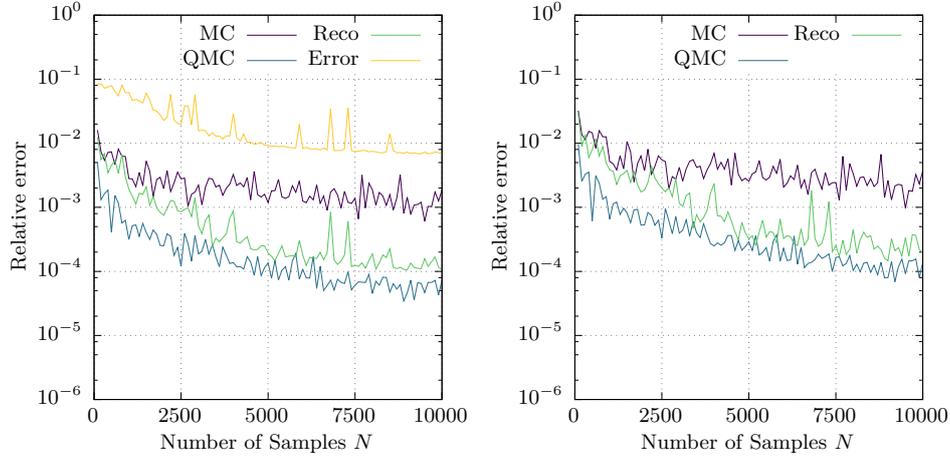

**Figure 6.** PDE setting (VI): TT reconstruction error of the *expectation* ("Reco") compared to MC and QMC simulations (left). Additionally the average relative error for a random set of 1000 samples *not* used for the reconstruction is shown ("Error"). TT reconstruction error of the *variance* compared to MC and QMC simulations (right).

Numerical experiments show that the pursued strategy provides a promising approach for solving high-dimensional PDEs and problems in UQ. A notable observation is that the achieved accuracy often is significantly better than what could be expected from the theoretical estimates. This raises the question if there are ways to improve the derivations for the problems under consideration. In particular, we have shown that for model classes in finite dimensional vector spaces the *generalization error* behaves like $N \in \mathcal{O}(\varepsilon^{-2} \ln \varepsilon)$. In fact, in the presence of noise in the input data, we should not expect major improvements. However, in the present setting, we are mainly interested in the noise-free case[2], since the samples are only used to approximate an integral. In a series of recent publications [13, 57] it was shown that in this case one can achieve the (almost) best approximation rate in $\ell^\infty$ with high probability for the approximation of problems like in this work. The authors considered linear model classes of (orthogonal) polynomials in $\mathbb{R}^d$ and the least squares loss function. Convergence is guaranteed by bounding the condition number of the Gram matrices. Using recent results of Chernóv type it can be shown that this holds for different types of polynomials high probability. A related analysis in the context of regression with sparse grids can be found in [32].

When compared to the works of Cucker et al., the analytical approach of Cohen et al. is in line with the fundamental theory of numerical methods for elliptic partial differential equation, manifested e.g. in Finite Elements [58] and projection methods [59]. Broadly speaking, convergence follows from *stability* together with *consistency*. Therefore, the approach of these authors is rooted mathematically in numerical analysis and based on a *stability* result, which in principle is an *inf-sup condition* [60], but limited to linear models.

We believe that this is an intriguing ansatz, which could also be adopted to non-linear model classes and empirical functionals, allowing for an explanation of our empirical observations. A way to extend the present results may be to formulate a corresponding RIP (restricted isometry property) condition for non-linear model classes due to which stability would follow immediately. This is e.g. common in the area of compressed sensing. From this perspective we conjecture that under the assumption of a RIP condition, the solution can be reconstructed exactly for a rather general set of model classes. This then allows for better estimates of the total error than what is achieved with the techniques used in this work. We defer an extended treatment of these ideas to forthcoming work. Nevertheless, the presented theory already provides a robust and versatile framework which guarantees convergence and is not restricted to linear models classes although the results appear to be too pessimistic for data with small noise.

---

[2]i.e. "noise" in the samples is only due the numerical approximation

We would like to point out that the condition $\ell(\Phi)(x) \leq C$ a.e. $\forall \Phi \in \mathcal{M}$ also asserts a certain kind of stability. However, in many cases it not easy to immediately show that this condition holds. In the case of RKHS, this estimate is guaranteed by the assumption $\|\Phi\|_\mathcal{K} \leq R$. If no regularization of a similar kind is invoked, it is not even clear for linear models that this condition holds. For example, this assumption can be violated in finite dimensional spaces if the dimension is larger than the number of samples.

In this article we neglect the optimization task, which otherwise raises many further question like the convergence to a global minimum and the complexity of the employed algorithm. Both aspects are of major importance for our PDE setting since the best approximation problem on tensor manifolds exhibits many local minima and is delicate with respect to computational complexity by its high-dimensional nature. In high-dimensional empirical risk minimization, variance reduced stochastic gradient methods [61–63] are established as a standard but an efficient adaptation to tensor networks is not straightforward. In the implementations it is also of interest to explore the possibility of importance sampling techniques (that may also be used to adaptively modify the sampling density) [14, 15] and of size adaptive sampling strategies [64, 65]. Moreover, recent experiences in the deep learning setting indicate that for highly non-convex models, stochastic gradient updates have an advantageous effect on the generalization error. This important issue is not understood yet and it is not clear how it impacts solving high-dimensional PDEs where the noise level is relatively low and the required accuracy is relatively high.

In order to render the presented approach more relevant for practical applications, adaptive strategies will have to be devised, which automatically steer the discretization parameters depending on the considered problem. Since the derived a priori bounds are quite pessimistic, they are inappropriate to yield an indication of the required sample number $N$. Moreover, in initial iterations, one does not need a very good approximation of the functional, which leaves some freedom to numerical methods.

To sketch some first ideas, recall that for given $\Phi \in \mathcal{M}$ it holds $\mathbb{E}[|\mathcal{J}(\Phi) - \mathcal{J}_N(\Phi)|] \leq \frac{C_1^4}{N^2}$. This provides a crude bound for the error of the approximation of the cost functional. Since optimizing $\mathcal{J}_N$ becomes more difficult with every sample and numerical (non-convex)

optimization is bound to errors, we argue that we can start with a relatively small set of samples and increase the sample size dynamically. This conjecture is confirmed in numerical experiments in Section 6. The reasoning is that if the optimization algorithm converges with order $\mathcal{O}(k^{-1})$ on $\mathcal{J}$, we want the error $|\mathcal{J} - \mathcal{J}_N|$ to be at least of the same order in some sense[3]. However, this bound for $\mathbb{E}[\mathcal{E}(\Phi)]$ cannot be used to choose $N$ since the constants are unknown. To still use this a posteriori bound, one could assume a minimization algorithm which converges with rate $k^{-1}$ (where $k$ is the iteration number). One can then choose a sequence $N_k$ s.t. the error in the approximation $\mathbb{E}[\mathcal{E}(\Phi)(N_k)] \in \mathcal{O}(k^{-1})$ vanishes as quickly as the error of the optimization.

Additionally, with the employed FE solution samples, an a posteriori error control of the approximation error can be achieved easily. If the approximation error is defined as the expectation of the energy error, an averaging of pathwise standard FE error estimators yields a sensible a.s. reliable error estimation for an adaptive mesh refinement procedure. The efficacy of this approach was presented in [66] for goal-oriented error control.

---

[3] such as with respect to the expectation as derived above